\documentclass[10pt,english]{extarticle}

\usepackage{caption}
\usepackage{enumerate}
\usepackage{t1enc}
\usepackage{babel}
\usepackage{amsthm}
\usepackage{amsmath}
\usepackage{amsfonts}
\usepackage{mathrsfs}
\usepackage[cp1250]{inputenc}
\usepackage{amsthm}
\usepackage{amsmath}
\usepackage{amssymb}
\usepackage{amsfonts}
\usepackage{indentfirst}
\usepackage{textcomp}
\usepackage{mathabx}
\usepackage{bbm}
\usepackage{tikz}
\usepackage{xcolor}
\usepackage{extsizes}
\usepackage[textwidth=18.5cm,textheight=23cm,left=3.5cm,voffset=-1cm,,footskip=1.5cm,headheight=0cm]{geometry}
\usepackage{fancyhdr}
\usepackage{yfonts}
\usepackage{float}
\usepackage{pgf}
\usepackage{pgfplots}
\usepackage{amsfonts,amsmath,amssymb,graphicx,array,url}
\usepackage{polski}
\usepackage[cp1250]{inputenc}
\pgfplotsset{compat=1.12}

\makeatletter
\def\namedlabel#1#2{\begingroup
    #2%
    \def\@currentlabel{#2}%
    \phantomsection\label{#1}\endgroup
}
\makeatother

\title{\LARGE \bf
Stability of long run functionals with respect to stationary Markov controls*
}

\author{Lukasz Stettner %
\thanks{*This work was supported by NCN grant 2020/37/B/ST1/00463}
\thanks{Lukasz Stettner is with the Institute of Mathematics Polish Academy of Sciences,
        Sniadeckich 8, 00-656 Warsaw, Poland,
        {\tt\small stettner@impan.pl}}%
}

\begin{document}

\maketitle
\thispagestyle{empty}
\pagestyle{empty}

\begin{abstract}

In the paper we study dependence of long run functionals and limit characteristics assuming that Borel measurable Markov controls converge pointwise. We consider two kinds of functionals: average cost per unit time and long run risk sensitive. We impose uniform ergodicity assumption, which is later is relaxed and suitable convergence of controlled transition probabilities.
\end{abstract}

\section{INTRODUCTION}
Assume that state space $E$ is a Polish space with Borel $\sigma$-filed ${\cal E}$. We have also a compact set of control parameters $U$ and a family ${\cal U}$ of Borel measurable mappings $u: E \mapsto U$ called later Markov controls. For each $u\in {\cal U}$ we are given a controlled Markov process $(X_t^u)$ with transition operator $P^{u(x)}(x,dy)$ for $x\in E$. We consider a natural pointwise convergence topology on ${\cal U}$, which means that $u_n\in {\cal U}$ converges to $u\in {\cal U}$ whenever $u_n(x)\to u(x)$ as $n\to \infty$ for each $x\in E$. In the significant part of the paper we shall assume that uniform ergodicity assumption is satisfied
\begin{equation}
(UE) \ \   \sup_{x,x'\in E} \sup_{a,a'\in U} \sup_{B\in {\cal E}} P^a(x,B)-P^{a'}(x',B):=\Delta <1
\end{equation}
where $P^a(x,\cdot)$ is a transition kernel with constant control $u\equiv a$. From section 5 of Chapter V in \cite{Doob} for any $u\in {\cal U}$ there is a measure $\pi^u\in {\cal P}(E)$ - the set of probability measures on $E$ such that for $x\in E$
\begin{equation}\label{eq0}
\sup_{B\in {\cal E}}|(P^u)^n(x,B)-\pi^u(B)|\leq \Delta^n
\end{equation}
where $(P^u)^n(x,\cdot)$ stands for $n$-th iteration of the transition kernel $(P^{u(x)}(x,\cdot))$. Clearly $\pi^u$ is a unique invariant measure for $(P^{u(x)}(x,\cdot))$, which means that $\int_E P^{u(x)}(x,B)\pi^u(dx)=\pi^u(B)$ for $B\in {\cal E}$. Assumption (UE) is restrictive in the case non compact state space since it requires that process almost immediately enters sufficiently large compact set. In the case of compact state space it more or less says that the process is mixing and it is satisfied when transition probability has a continuous positive density function. To have property (\ref{eq0}) is suffices to have (UE) satisfied for a certain iteration of transition probabilities.   
Assume that $u_n,u\in {\cal U}$ and $u_n\to u$. We want to find sufficient conditions for the following continuity results

\noindent
{\bf Problem 1.} $\sup_{B\in {\cal E}}|\pi^{u_n}(B)-\pi^u(B)|\to 0$ as $n\to \infty$,

\noindent
{\bf Problem 2.} $J_x(u_n)\to J_x(u)$, as $n\to \infty$, where
\begin{equation}
J_x(u)=\liminf_{n\to \infty}{1\over n} E_x^u\left\{\sum_{i=0}^{n-1} c(X_i^u,u(X_i^u))\right\}
\end{equation}
 for a bounded measurable function $c: E\times U\mapsto R$ continuous with respect to the second (control) parameter.

\noindent
 Let for $B\in {\cal B}(E)$ and $\alpha \in (-\infty,+\infty)\setminus\left\{0\right\}$
 \begin{equation}
 \lambda_x^{u,\alpha}(B)=\liminf_{n\to \infty} {1\over \alpha n} \ln E_x^u\left\{e^{\alpha \sum_{i=0}^{n-1} 1_B(X_i^u)}\right\}
 \end{equation}

 \noindent
 {\bf Problem 3.} $\sup_{B\in {\cal E}} |\lambda_x^{u_n,\alpha}(B)-\lambda_x^{u,\alpha}(B)|\to 0$ as $n\to \infty$,
for $\alpha\in (-\infty,+\infty)\setminus\left\{0\right\}$,

 \noindent
 {\bf Problem 4.} $I_x^\alpha(u_n)\to I_x^\alpha (u)$, as $n\to \infty$, where
 \begin{equation}
 I_x^\alpha(u)=\liminf_{n\to \infty} {1\over n} \ln E_x^u\left\{e^{\alpha \sum_{i=0}^{n-1} c(X_i^u,u(X_i^u))}\right\},
\end{equation}

\noindent
{\bf Problem 5.} $\lim_{\alpha \to 0} \sup_{B\in {\cal E}} | \lambda_x^{u,\alpha}(B)- \pi^u(B)|= 0$

\noindent
{\bf Problem 6.} for $\alpha_n\to 0$, ${1\over \alpha_n}I_x^{\alpha_n}(u_n)\to J_x(u)$, as $n\to \infty$.

Problem 2 corresponds to average reward per unit time. Under usual ergodic assumptions (see \cite{HLL}) the value of $J_x(u)$ is constant not depending on $x$. It is in fact an integral with respect to an invariant measure. Similarly one can expect that value $\lambda_x^{u,\alpha}$ does not depend on $x$. Problem 4 has an interpretation as risk sensitive control (see \cite{DiMS1} for motivations) with risk factor $\alpha$, which maybe positive or negative. It measures expected value of the reward functional with further moments with weights determined by the risk factor $\alpha$.  To solve problems 1-4, in particular 3-4 we have to impose nice ergodic structure of the controlled Markov process. Uniform egodicity (UE) introduced in \cite{Doob} is a strong assumption (see also \cite{MT}). We relax it in the case of the problem 2 in section IV. What is important for us is that we do not require Feller property (which usually is not satisfied for discrete time controlled Markov processes) of the controlled Markov processes and therefore we work on approximation of the limit measures in variation norm. A typical result for controlled Feller Markov processes can be formulated as follows

\noindent
{\bf Lemma 1.} {\it  Assume that for $f\in C(E)$  the space of continuous bounded functions on $E$, the mapping
$E\times U \ni (x,a)\mapsto \int_E f(y) P^a(x,dy)$ is continuous, for each continuous $u:E\mapsto U$ there is a unique invariant measure $\pi^u$. Assume furthermore that for continuous functions $u_n$ and $u$ we have $u_n(x)\to u(x)$ as $n\to \infty$ and the measures $\pi^{u_n}$ are tight. Then  $\pi^{u_n}(f)\to \pi^u(f)$ for any $f\in C(E)$.}

\noindent
{\bf Proof.} We have that
\begin{equation}\label{row}
\pi^{u_n}(f)=\int_E\int_E f(y)P^{u_n(x)}(x,dy)\pi^{u_n}(dx)
\end{equation}
Since $\pi^{{u_n}}$ are tight by Prokhorov theorem (see \cite{Bil} thm. 6.1) there is a subsequence, for simplicity still denoted by $\pi^{u_n}$ and a measure $\pi$ such that $\pi^{{u_n}}(f)\to \pi(f)$ for $f\in C(E)$. Moreover for each $\epsilon>0$ there is a compact set $K$ such that $\pi^{u_n}(K)\geq 1-\epsilon$ and $\pi(K)\geq 1-\epsilon$ .
Now
\begin{eqnarray}
&& |\int_E\int_E f(y)P^{u_n(x)}(x,dy)\pi^{u_n}(dx)- \int_E\int_E f(y)P^{u(x)}(x,dy)\pi(dx)|\leq \nonumber \\
&& 2\|f\|\epsilon + \sup_{x\in K}|\int_E f(y)P^{u_n(x)}(x,dy)-\int_E f(y)P^{u(x)}(x,dy)| 
\end{eqnarray}
so that letting $n\to \infty$ in (\ref{row}) we obtain
\begin{equation}\label{ro}
\pi(f)=\int_E\int_E f(y)P^{u(x)}(x,dy)\pi(dx)
\end{equation}
for $f\in C(E)$. This means that $\pi$ is an invariant measure for $P^{u(x)}(x,\cdot)$. By uniqueness we have that $\pi=\pi^u$. Since for any other convergent subsequence $\pi^{u_n}$ we also have $\pi^u$ as a weak limit measure, we finally have that $\pi^{u_n}(f)\to \pi^u(f)$ for $f\in C(E)$ and $n\to \infty$.

As one can see the restriction to continuous controls (required by the weak convergence techniques see e.g. \cite{KD}) is important. In this paper we try omit this using uniform ergodicity of Markov or embedded Markov process. 

Problems 1-6 arise naturally when we want to solve stochastic control problems. Usually we are not able to solve suitable Bellman equation explicitly. By general theory (see \cite{HLL} in the average reward per unit time problem and \cite{DiMS1} for long run risk sensitive problem) we can show that optimal control is Markov i.e. it is a function of the current value of the state process. Such functions are usually only Borel measurable and a natural question is how to approximate them. In the paper we address this problem for long run functionals. Our main result is that pointwise approximation of the control function is stable with respect to limit (invariant) measures and this way we approximate the value function as well. Therefore we seek nearly optimal control within a class of piecewise constant controls which approximate it an optimal  control. Such class of controls appears naturally when we consider state space discretization. Our continuity results also  justify commonly used procedure,  when instead of solving suitable (discretized) Bellman equations we consider Monte Carlo simulations and improving piecewise constant controls we want to obtain a reasonable approximation of nearly optimal control. Potential limit of such procedure using finer and finer discretization steps leads us (thanks to shown below continuity of functionals) to an optimal Markov control. 
Problems 5 and 6 concern asymptotics of risk sensitive functionals with respect for sufficiently small value of risk parameter. For this purpose we adopt Hoeffding's lemma (see Lemma 2.6 of \cite{Ma}). Average reward per unit can be approximated by risk sensitive problems with small risk and on the other hand average reward per unit time can be considered as an approximation of risk sensitive problem with sufficiently small risk factor. What is important we solve problem 5 and 6 without requiring continuity of controls, which was crucial when we used large deviations of empirical measures result (see \cite{DiMS1} or  \cite{LS2}) and this is again another method to construct nearly optimal Markov controls.

\section{Average reward per unit time under (\ref{eq0})}

We start with the following important

\noindent
{\bf Proposition 1.}
{\it Assume that $\mu_n,\mu \in {\cal P}(E)$ are such that $\sup_{B\in {\cal E}}|\mu_n(B)-\mu(B)|\to 0$, as $n\to \infty$ and $f_n:E\mapsto R $ is a bounded sequence of bounded Borel measurable functions such that $f_n(x)\to f(x)$, as $n\to \infty$ for $x\in E$. Then
$\mu_n(f_n)\to \mu(f)$, as $n\to \infty$.}

\noindent{\bf Proof.} Without loss of generality we may assume that $0\leq f_n\leq 1$. For a given $\epsilon>0$ let $N$ be a positive integer such that ${1\over N}\leq \epsilon$. Define set
$ F_n^m=\left\{x: {m-1 \over N}\leq f_n(x)<{m\over N}\right\} $ for $m<N$ and $F_n^N=\left\{x: {N-1\over N} \leq f_n(x)\leq 1\right\}$. Then taking into account that $\sup_{x\in F_n^m}f_n(x)-\inf_{x\in F_n^m}f_n(x)\leq {1\over N}$ we have
\begin{equation}
\sum_{m=1}^N |f_n(y)\mu_n(dy) - {m-1\over N}\mu_n(F_n^m)|\leq {1\over N}
\end{equation}
and the same holds when we replace $\mu_n$ by $\mu$ we have 
\begin{eqnarray}\label{in1}
&&|\mu_n(f_n)-\mu(f_n)|\leq |\sum_{m=1}^N \int_{F_n^m} f_n(y)(\mu_n(dy)-\mu(dy))|\leq {2\over N} \nonumber \\
&& + \sum_{m=1}^N {m-1 \over N} \sup_{B\in {\cal E}}|\mu_n(B)-\mu(B)|\leq 2\epsilon + {N-1\over 2} \sup_{B\in {\cal E}}|\mu_n(B)-\mu(B)|\to 2 \epsilon
\end{eqnarray}
as $n\to \infty$. Now
\begin{equation}\label{in2}
|\mu_n(f_n)-\mu(f)|\leq |\mu_n(f_n)-\mu(f_n)|+|\mu(f_n)-\mu(f)|\to 2 \epsilon
\end{equation}
as $n\to \infty$ and since $\epsilon$ could be arbitrarily small we complete the proof.

The next Proposition uses several arguments of the proof of Proposition 1.

\noindent
{\bf Proposition 2.} {\it  Assume that $u_n\in {\cal U}$ converges to $u\in {\cal U}$ and
\begin{equation}\label{m1}
\sup_{B\in {\cal E}} |P^{u_n(x)}(x,B)-P^{u(x)}(x,B)|\to 0
\end{equation}
as $n\to \infty$ for $x\in E$. Then for each positive integer $k$
\begin{equation}\label{gm}
\sup_{B\in {\cal E}} |(P^{u_n})^k(x,B)-(P^{u})^k(x,B)|\to 0
\end{equation}
as $n\to \infty$ for $x\in E$.}

\noindent
{\bf Proof.} We use induction. For $k=1$ (\ref{gm}) follows directly from (\ref{m1}). Assume that (\ref{gm}) is satisfied for $k$. Then for $k+1$ we have
\begin{equation}
(P^{u_n})^{k+1}(x,B)=\int_E (P^{u_n})^k(y,B)P^{u_n(x)}(x,dy)
\end{equation}
Let $f_n(y,B)= (P^{u_n})^k(y,B)$, $f(y,B)=(P^{u})^k(y,B)$, $\mu_n(\cdot)=P^{u_n(x)}(x,\cdot)$ and $\mu(\cdot)=P^{u(x)}(x,\cdot)$.  By induction hypothesis $\sup_{B\in {\cal E}} |f_n(y,B)-f(y,B)|\to 0$ and
$\sup_{B\in {\cal E}} |\mu_n(B)-\mu(B)|\to 0$ as $n\to \infty$. Moreover $f_n(y,B)$ takes values in the interval $[0,1]$. We follow now the proof of Proposition 1. For $\epsilon>0$ take positive integer $N$ such that ${1\over N}\leq \epsilon$. Define sets $F_n^k$ and $F^k$ as in the proof of Proposition 1 (depending now on the set $B$). Then by (\ref{in1})
\begin{equation}
\sup_{B\in {\cal E}}|\mu_n(f_n(B))-\mu(f_n(B))|\leq 2\epsilon + {N-1\over 2} \sup_{A\in {\cal E}}|\mu_n(A)-\mu(A)|\to 2\epsilon
\end{equation}
as $n\to \infty$.
Now as in (\ref{in2})
\begin{equation}
\sup_{B\in {\cal E}}|\mu_n(f_n(B))-\mu(f(B))|\leq \sup_{B\in {\cal E}}|\mu_n(f_n(B))-\mu(f_n(B))|+\sup_{B\in {\cal E}}|\mu(f_n(B))-\mu(f(B))|\to 2\epsilon
\end{equation}
as $n\to \infty$, which completes the proof.

We can now solve problem 1 and problem 2 under (\ref{eq0})

\noindent{\bf Theorem 1.} {\it Assume that (\ref{eq0}) is satisfied, $u_n\in {\cal U}$ converges to $u \in {\cal U}$ and we have (\ref{m1}). Then
\begin{equation}\label{prob1}
\sup_{B\in {\cal E}} |\pi^{u_n}(B)-\pi^u(B)|\to 0
\end{equation}
 and
 \begin{equation}\label{prob2}
J_x(u_n)=\int_E c(y,u_n(y))\pi^{u_n}(dy)\to J_x(u)=\int_E c(y,u(y))\pi^{u}(dy)
 \end{equation}
 as $n\to \infty$.}

 \noindent{\bf Proof.} From (\ref{eq0}) we have that
 \begin{equation}
 \sup_{B\in {\cal E}}|(P^{u_n})^k(x,B)-\pi^{u_n}(B)|\leq \Delta^k
 \end{equation}
 and
 \begin{equation}
 \sup_{B\in {\cal E}}|(P^u)^k(x,B)-\pi^{u}(B)|\leq \Delta^k.
 \end{equation}
 Therefore
 \begin{equation}
 \sup_{B\in {\cal E}}|\pi^{u_n}(B)- \pi^{u}(B)|\leq \sup_{B\in {\cal E}}|(P^u_n)^k(x,B)-(P^u)^k(x,B)|+2 \Delta^k.
 \end{equation}
 By (\ref{gm}) for fixed $k$ letting $n\to \infty$ we obtain that
 $\sup_{B\in {\cal E}}|\pi^{u_n}(B)- \pi^{u}(B)|\leq \Delta^k$, and letting $k\to \infty$ we have the claim. Having (\ref{prob1}) immediately from Proposition 1 we obtain (\ref{prob2}). The proof is completed.

\noindent{\bf Remark 1.} {\it We can relax assumption (\ref{eq0}) imposing that 
\begin{equation}
\sup_{x,x'\in E} \sup_{a,a'\in U} {\int_E V(y)|P^a(x,dy)-P^{a'}(x',dy)| \over V(x)+V(x')}:=\Delta<1
\end{equation}
for a Borel measurable function $V:E\mapsto [1,\infty)$ such that for some $x\in E$ we have $\sup_{a\in U} P^aV(x)<\infty$. Then by section 10.2c of \cite{HLL} there is $R>0$ such that for any $u\in {\cal U}$, $x\in E$  
\begin{equation}
\sup_{B\in {\cal E}} |(P^u)^k(x,B)-\pi^u(B)|\leq \Delta^k R V(x)
\end{equation}
and we can get the same claim as in Theorem 1.} 

\section{Risk sensitive control}
In this section we shall need the following assumption

\noindent
(ME) \ \  there is an integer $m\geq 1$ and a constant $K$ such that
\begin{equation*}
\sup_{x,x'\in E} \sup_{u\in {\cal U}} \sup_{B\in {\cal E}} {(P^u)^m(x,B) \over (P^u)^m(x',B)}:=K<\infty 
\end{equation*}

Assumption (ME) is quite restrictive in the case of noncompact state spaces. It says that $m$-th iterations of  transition probabilities are equivalent with densities bounded from above and separated from $0$. When the state space is compact and densities of the $m$-th iterations of transition probabilities are continuous and positive, assumption (ME) is clearly satisfied.    
Let $B(E)$ be the set of bounded Borel measurable functions on $E$ with supremum norm. For $g\in B(E)$ define so called span norm  $\|g\|_{sp}=\sup_{x\in E} g(x)- \inf_{x'\in E} g(x')$. Let for $u\in {\cal U}$ and $f, g\in B(E)$, and $\alpha\in (-\infty, + \infty)\setminus \left\{0\right\}$
\begin{equation}
\Psi^{u,\alpha} g(x)= f(x) + {1 \over \alpha}\ln \int_E e^{\alpha g(y)}P^{u(x)}(x,dy)
\end{equation}
We have

\noindent
{\bf Theorem 2.} {\it  Under (UE) operator $\Psi^{u,\alpha}$ is a local contraction in the span norm in the space $B(E)$ for $u\in {\cal U}$, i.e. there is a function $\gamma_\alpha: (0,\infty) \mapsto [0,1)$ such that whenever for $g_1,g_2 \in B(E)$ we have $\|g_1\|_{sp}\leq M$ and $\|g_2\|_{sp}\leq M$ then
\begin{equation}
\|\Psi^{u,\alpha} g_1 - \Psi^{u,\alpha} g_2\|_{sp}\leq \gamma_\alpha(M)\|g_1-g_2\|_{sp}.
\end{equation}
Furthermore additionally under (ME) operator  $\Psi^{u,\alpha}$ is a global contraction in $B(E)$ with the span norm and the space $B(E)$ is transformed by $\Psi^{u,\alpha}$ into the subspace of $B(E)$ with the span norm less that $\tilde{K}$, where $\tilde{K}$ does no depend on $u\in {\cal U}$ but depends on $K$ from (ME).}

\noindent
{\bf Proof.}  Local contractivity follows from Theorem 3, Corollary 4 and 5 in \cite{LS2}. We give here few hints. Using dual representation of the operator $\Psi$ (see Proposition 1.42 of \cite{DE}) we have that 
\begin{equation}\label{er1}
\Psi^{u,\alpha} g(x)= f(x) + \inf_{\nu\in {\cal P}(E)}\left[\int_E g(y)\nu(dy)-{1\over \alpha}H(\nu,P^{u(x)}(x,\cdot))\right]
\end{equation}
with $H(\nu_1,\nu_2):=\int_E \ln({d\nu_1 \over d \nu_2})d\nu_1$ when $\nu_1$ is absolutely continuous with respect to $\nu_2$, and is equal to $+\infty$ otherwise. 
When $\alpha<0$ with infimum attained for the measure $\nu_{x,u,\alpha g}(B)={\int_Be^{\alpha g(y)}P^{u(x)}(x,dy)\over 
\int_Ee^{\alpha g(y)}P^{u(x)}(x,dy)}$. Infimum is replaced by maximum when $\alpha>0$ and we have then the same maximizing measure. 
For $g_1,g_2\in B(E)$ and $x_1,x_2\in E$ using (\ref{er1}) we obtain
\begin{equation}
\Psi^{u,\alpha}g_1(x_1)-\Psi^{u,\alpha}g_2(x_1)-\Psi^{u,\alpha}g_1(x_2)+ \Psi^{u,\alpha}g_1(x_2) \leq\|g_1-g_2\|_{sp} \sup_B (\nu_1-\nu_2)(B)
\end{equation}
Now using the form of measures $\nu_{x_i,u,\alpha g_i}$ for $i=1,2$, by assumption (UE) we obtain Lipschitz constant $\gamma_\alpha(M)$.  
Global contraction then follows from Remark 4 and Proposition 6 in \cite{LS1}.

We now have

\noindent
{\bf Corollary 1.} {\it For $u\in {\cal U}$ and fixed $f\in B(E)$ there is a constant $\lambda^{u,\alpha}(f)$ and a function $w_f^{u,\alpha}\in B(E)$ such that for $x\in E$ we have
\begin{equation}\label{riskbel}
\alpha w_f^{u,\alpha}(x)=\alpha (f(x)-\lambda^{u,\alpha}(f)) + \ln \int_E e^{\alpha w_f^{u,\alpha}(y)}P^{u(x)}(x,dy)
\end{equation}
Moreover $\|w_f^{u,\alpha}\|_{sp}\leq \tilde{K}$, where $\tilde{K}$ depends on $m$ and $K$ from (ME) and the function $\gamma_\alpha$.}

\noindent
{\bf Proof.} By Theorem 2 there is a fixed point $w_f^{u,\alpha}$ of the operator $\Psi^{u,\alpha}$ i.e. $\|\Psi^{u,\alpha} w_f^{u,\alpha}-w_f^{u,\alpha}\|_{sp}=0$. Therefore there is a constant $\lambda^{u,\alpha}(f)$ such that $\Psi^{u,\alpha} w_f^{u,\alpha}(x)-\lambda_f^{u,\alpha} = w_f^{u,\alpha}$, which completes the proof.

\noindent
Iterating (\ref{riskbel}) and noticing that $\|w_f^{u,\alpha}\|_{sp}\leq \tilde{K}$ we obtain

\noindent{\bf Corollary 2.} {\it  For a positive integer $k$ we have for $x\in E$, $u\in {\cal U}$, $f\in B(E)$
\begin{equation}
\alpha w_f^{u,\alpha}(x)= - k \alpha \lambda^{u,\alpha}(f)+ \ln E_x^u\left\{e^{\alpha(\sum_{i=0}^{k-1} f(X_i^u)+ w_f^{u,\alpha}(X_k))}\right\}
\end{equation}
and for $u_n,u\in {\cal U}$
\begin{equation}\label{est1}
k|\lambda^{u_n,\alpha}(f)-\lambda^{u,\alpha}(f)|\leq |{1\over \alpha}\ln {E_x^{u_n}\left\{e^{\alpha \sum_{i=0}^{k-1}f(X_i^{u_n})}\right\}\over
E_x^{u}\left\{e^{\alpha \sum_{i=0}^{k-1}f(X_i^{u})}\right\}}|+ 2\tilde{K}.
\end{equation}}

Therefore we obtain

\noindent
{\bf Proposition 3.} {\it Under (\ref{m1}) we have that when $u_n\to u$ and $u_n,u\in {\cal U}$ then for positive integer $k>1$ we have
\begin{equation}\label{riskm}
\sup_{B\in {\cal E}} |\ln {E_x^{u_n}\left\{e^{\alpha \sum_{i=0}^{k-1}1_B(X_i^{u_n})}\right\}\over
E_x^{u}\left\{e^{\alpha \sum_{i=0}^{k-1}1_B(X_i^{u})}\right\}}|\to 0
\end{equation}}

\noindent
{\bf Proof.} Notice first that it is sufficient to show that for any $r>0$ and sequence $f_i, f_i^n\in B(E)$ such that $\|f_i^n\|=\sup_{y\in E} |f_i^n(y)|\leq r$, $f_i^n(y)\to f_i(y)$ as $n\to \infty$
\begin{equation}\label{indmul}
\sup_{f_i^n, f_i \in B(E), \|f_i^n\|, \|f^n\|\leq r} |E_x^{u_n}\left\{e^{\sum_{i=0}^{k-1}f_i^n(X_i^{u_n})}\right\} -  E_x^{u}\left\{e^{\sum_{i=0}^{k-1}f_i(X_i^{u})}\right\}|\to 0
\end{equation}
as $n\to \infty$.
We use induction. For $k=1$ (\ref{indmul}) is clearly satisfied. Assume (\ref{indmul}) for $k$. Then
\begin{equation}
E_x^{u_n}\left\{e^{\sum_{i=0}^{k}f_i^n(X_i^{u_n})}\right\}= E_x^{u_n}\left\{e^{\sum_{i=0}^{k-2}f_i^n(X_i^{u_n})+\tilde{f}_{k-1}^n}\right\}
\end{equation}
with $\tilde{f}_{k-1}^n(y)=f_{k-1}^n(y) + \ln \int_E e^{f_k^n(z)}P^{u_n(y)}(y,dz)$.
Similarly
\begin{equation}
E_x^{u}\left\{e^{\sum_{i=0}^{k}f_i(X_i^{u})}\right\}=E_x^{u}\left\{e^{\sum_{i=0}^{k-2}f_i(X_i^{u})+\tilde{f}_{k-1}}\right\}
\end{equation}
with $\tilde{f}_{k-1}(y)=f_{k-1}(y) + \ln \int_E e^{f_k(z)}P^{u(y)}(y,dz)$.
Therefore using (\ref{m1}) we have that $\tilde{f}_{k-1}^n(y) \to \tilde{f}_{k-1}(y)$, as $n\to \infty$ for $y\in E$ and consequently by induction hypothesis (\ref{indmul}) holds for $k+1$. Therefore (\ref{riskm}) is satisfied.

\noindent
From Corollary 2 taking into account Proposition 3 we obtain

\noindent
{\bf Theorem 3.} {\it Under (UE), (ME) and (\ref{m1}) for ${\cal U}\ni u_n\to u\in {\cal U}$ and $\alpha \in (-\infty, +\infty) \setminus\left\{0\right\}$ we have
\begin{equation}\label{p3}
\sup_{B\in {\cal B(E)}} |\lambda_x^{u_n,\alpha}(B)-\lambda_x^{u,\alpha}(B)|\to 0
\end{equation}
 as $n\to \infty$. Furthermore
$I_x^\alpha(u_n) \to I_x^\alpha(u)$ as $n\to \infty$.}

\noindent
{\bf Proof.} From (\ref{est1}) using (\ref{riskm}) we obtain (\ref{p3}). From (\ref{riskm}) and Proposition 1 we have that
\begin{equation}
E_x^{u_n}\left\{e^{\alpha \sum_{i=0}^{k-1}c(X_i^{u_n},u_n(X_i^{u_n}))}\right\} \to E_x^{u}\left\{e^{\alpha \sum_{i=0}^{k-1}c(X_i^{u},u(X_i^{u}))}\right\}
\end{equation}
as $n\to \infty$. Therefore from (\ref{est1}) we obtain the convergence of risk sensitive functionals, which completes the proof.

\noindent
{\bf Remark 2.} {\it Assumption (UE) plays an important role to study risk sensitive Bellman equation. We can try to relax it using some splitting technics from the paper \cite{DiMS2}. We then require a number of additional assumptions and therefore such results are far away from the scope of this paper.
Assumption (ME) can be replaced by requiring small risk $|\alpha|$ as was studied in the papers \cite{DiMS3} or \cite{PS}.
Using assumption (UE) we are looking for a bounded solution to (\ref{riskbel}). We can use other technics based on Krein Rutman theorem (see \cite{LS3} and \cite{A}) or suitable Lyapunov conditions (see \cite{Ch}) and work with unbounded solutions. In such case and analog of Theorem 3 would require more assumptions.    }

\section{Average reward per unit time without (UE)}
Assume we are give two concentric balls $R$ and $R_1$, $R\subset R_1$. Define
$D_B=\inf\left\{s\geq 0: X_s\in B\right\}$ and $\tau_B=D_{R_1^c}+D_B\circ \theta_{D_{R_1^c}}$. Stopping time $\tau_B$ is the first time when process enters the set $B$ after entering the complement of the ball $R_1$. We shall assume that

\noindent
\begin{eqnarray*}
(EC) \ \ \ \ \ \ &&  \sup_{x\in R} \sup_{u\in {\cal U}} E_x^u\left\{(\tau_R)^2\right\}<\infty. \ \ \ \  \  \ \ \  \ \ \ \ \ \ \ \ \ \ \ \ \
\end{eqnarray*}

For $x\in R$ and $B\in {\cal B}(R)$ - the family of Borel subsets of $R$ let
\begin{equation}
\Pi^u(x, B):=P_x^u\left\{X_{\tau_R}^u\in B\right\}.
\end{equation}
We assume that embedded Markov process with transition operator $\Pi^u$ is uniformly ergodic i.e.

\noindent
\begin{eqnarray*}
(ER)  \ \ \ \ && \sup_{x, x'\in R} \sup_{u\in {\cal U}}\sup_{B\in {\cal B}(R)} |\Pi^u(x,B)-\Pi^u(x',B)| :=\Delta_R<1
\end{eqnarray*}

From section 5 of Chapter V in \cite{Doob} for any $u\in {\cal U}$ there is a unique invariant measure $\mu^u\in {\cal P}(R)$  for the transition operator $\Pi^u$ such that for positive integer $n$ we have
\begin{equation}\label{oszac}
\sup_{B\in {\cal B}(R)} |(\Pi^u)^n(x,B)-\mu^u(B)|\leq (\Delta_R)^n.
\end{equation}
Assume that

\noindent
\begin{eqnarray*}
(PR) \ \ \ \ \ \ \  \ \ \ \ \ \ \ \sup_{x\in E} \sup_{u\in {\cal U}} E_x^u\left\{D_R\right\}<\infty.\ \ \ \ \ \ \ \ \ \ \ \ \ \ \ \ \
\end{eqnarray*}

Then following Chapter 3 or \cite{Kh} there is a unique invariant measure $\pi^u$ for $(X^u_t)$ and it is of the form for $B\in {\cal E}$
\begin{equation}
\pi^u(B)={\int_R E_x^u\left\{\sum_{i=0}^{\tau_R-1}1_B(X_i^u)\right\}\mu^u(dx) \over
\int_R E_x^u\left\{\tau_R\right\}\mu^u(dx)}.
\end{equation}
We shall need the following

\noindent
{\bf Proposition 4.} {\it Under (EC), (\ref{m1}), for $u_n,u\in {\cal U}$, $u_n\to u$ we have
\begin{equation}\label{mo}
\sup_{B\in {\cal B}(R)}|\Pi^{u_n}(x,B)-\Pi^u(x,B)|\to 0
\end{equation}
as $n\to \infty$.}

\noindent
{\bf Proof.} Notice first that for $x\in R$ and $B\in {\cal B(R)}$ we have
\begin{equation}\label{mom}
|\Pi^{u_n}(x,B)-P_x^{u_n}\left\{X_{\tau_R}^{u_n}\in B, \tau_R<N\right\}|\leq P_x^{u_n}\left\{\tau_R\geq N\right\}\leq {1\over N} E_x^{u_n}\left\{\tau_R\right\}.
\end{equation}
Now for $x\in R$, $2\leq N-k$ and $j\leq N-k-1$ we have
\begin{eqnarray}
&& P_x^{u_n}\left\{X_{\tau_R}^{u_n}\in B, D_{R_1^c}=j, \tau_R=N-k\right\}= \int_{R_1}\int_{R_1} \ldots \int_{R_1^c}\ldots \int_{R^c} \int_{R^c}\ldots \int_{R^c} P^{u_n(y_{N-k-1})}(y_{N-k-1},B) \nonumber \\
&& P^{u_n(y_{N-k-2})}(y_{N-k-2},dy_{N-k-1})\ldots P^{u_n(y_{j})}(y_{j},dy_{j+1}) P^{u_n(y_{j-1})}(y_{j-1},dy_j)\ldots P^{u_n(x)}(x,dy_1)
\end{eqnarray}
Under (\ref{m1}) we can show (detailed proof can be shown using induction as in the proof of Proposition 2) that
\begin{equation}
\sup_{B\in {\cal{B}}(R)} |P_x^{u_n}\left\{X_{\tau_R}^{u_n}\in B, D_{R_1^c}=j, \tau_R=N-k\right\}- P_x^{u}\left\{X_{\tau_R}^{u}\in B, D_{R_1^c}=j, \tau_R=N-k\right\}|\to 0
\end{equation}
as $n\to \infty$. Taking into account (\ref{mom}) we obtain (\ref{mo}).

In analogy to Proposition 2. we now have

\noindent
{\bf Corollary 3.} {\it Under (EC), (\ref{m1}), for $u_n,u\in {\cal U}$, $u_n\to u$ we have for positive integer $k$
\begin{equation}\label{Gm}
\sup_{B\in {\cal B}(R)} |(\Pi^{u_n})^k(x,B)-(\Pi^{u})^k(x,B)|\to 0
\end{equation}
as $n\to \infty$ for $x\in R$.}

As in the proof of Theorem 1 we now immediately have from (\ref{oszac})

\noindent
{\bf Corollary 4.} {\it Under (EC), (\ref{m1}), for $u_n,u\in {\cal U}$, $u_n\to u$ we have
\begin{equation}
\sup_{B\in {\cal B}(R)}|\mu^{u_n}(B)-\mu^u(B)|\to 0
\end{equation}
as $n\to \infty$.}

\noindent
{\bf Proposition 5.}  {\it  Under (EC), (\ref{m1}) when $u_n \to u$ for $u_n,u\in {\cal U}$ and bounded sequence $f_n\in B(E)$ is such that $f_n(x)\to f(x)$ for $n\to \infty$ then
\begin{equation}\label{os1}
E_x^{u_n}\left\{\sum_{i=0}^{\tau_R-1}f_n(X_i^{u_n})\right\} \to E_x^u\left\{\sum_{i=0}^{\tau_R-1}f(X_i^u)\right\}
\end{equation}
as $n\to \infty$ for $x\in E$.}

\noindent
{\bf Proof.} By (EC) it is sufficient to show that for positive integer $N$ we have
\begin{equation}
E_x^{u_n}\left\{\sum_{i=0}^{\tau_R-1}f_n(X_i^{u_n})1_{\tau_R<N}\right\} \to E_x^u\left\{\sum_{i=0}^{\tau_R-1}f(X_i^u)1_{\tau_R<N}\right\}.
\end{equation}
Now, it suffices to show that for each $1\leq j\leq N-1$
\begin{equation}\label{last}
E_x^{u_n}\left\{\sum_{i=0}^{\tau_R-1}f_n(X_i^{u_n})1_{D_{R_1^c}=j} 1_{\tau_R<N}\right\} \to  E_x^u\left\{\sum_{i=0}^{\tau_R-1}f(X_i^u)1_{D_{R_1^c}=j}1_{\tau_R<N}\right\}.
\end{equation}
The convergence (\ref{last}) can be shown by induction in a similar way as in the proof of Proposition 2 using Proposition 4. Therefore we have (\ref{os1}).

We can now summarize Proposition 5, Proposition 4 using Proposition 1.

\noindent
{\bf Theorem 4.}
{\it Under (EC), (ER), (PR) and (\ref{m1}) we have
\begin{equation}\label{r1}
sup_{B\in {\cal E}} |\pi^{u_n}(B)-\pi^u(B)|\to 0
\end{equation}
and
\begin{equation}\label{r2}
J_x(u_n) \to J_x(u)=\int_E c(x,u(x))\pi^u(dx)
\end{equation}
as $n\to \infty$.}

\noindent
{\bf Proof.} By Proposition 5 and Proposition 1 we obtain (\ref{r1}). Now, using (\ref{r1}), also (PR) and convergence of $c(x,u_n(x))\to c(x,u(x))$ we obtain (\ref{r2}). The proof is therefore completed.

\section{Risk sensitive asymptotics}

We shall now compare average reward per unit time functional with risk sensitive functional with small risk solving problems 5 and 6. We start with the following powerful Lemma

\noindent
{\bf Lemma 2.} {\it For a random variable $X$ taking values in the interval $[a,b]$ and any real $\alpha$  we have
\begin{equation}\label{Hoe}
0\leq \ln E\left\{e^{\alpha X}\right\}-\alpha E\left[X\right]\leq {(b-a)^2 \over 8} \alpha^2.
\end{equation}}

\noindent
{\bf Proof.} It follows from   Hoeffding's lemma (see Lemma 2.6 of \cite{Ma})  that
\begin{equation}
\ln E\left\{e^{\alpha (X-E\left[X\right])}\right\}\leq {(b-a)^2 \over 8} \alpha^2.
\end{equation}
Since the mapping $(-\infty,+\infty)\ni \alpha \mapsto {1 \over \alpha}\ln  E\left\{e^{\alpha X}\right\}$ with value $E[X]$ for $\alpha=0$ is increasing we therefore have that $\ln E\left\{e^{\alpha X}\right\}-\alpha E\left[X\right]\geq 0$, which completes the proof.

The following theorem solves problem 5

\noindent
{\bf Theorem 5.} {\it Under (UE) and (ME)  we have that for $x\in E$
\begin{equation}
\sup_{u\in {\cal U}} \sup_{B\in {\cal E}} |\lambda_x^{u,\alpha}(B)-\pi^u(B)|\to 0
\end{equation}
as $\alpha \to 0$.}

\noindent
{\bf Proof.} By Theorem 2 for such $\alpha$ and $f\in B(E)$ such that $0\leq f \leq 1$ there is a function $w_f^{u,\alpha}$ and constant $\lambda^{u,\alpha}(f)$ satisfying the equation (\ref{riskbel}). Furthermore analysis of the proof of Theorem 5 of \cite{LS1} shows that $\sup_{|\alpha|\leq \kappa}\gamma_\alpha(M)<1$ for $\kappa >0$ and any $M>0$. Consequently $\sup_{|\alpha|\leq \kappa}w_f^{u,\alpha}\leq D$ and $D$ does not depend on $u\in {\cal U}$ and depends only on the upper bound of the span norm of $f$. Therefore
\begin{eqnarray}\label{ossa}
&& |\lambda_x^{u,\alpha}(f)-\pi^u(f)|\leq |\lambda_x^{u,\alpha}(f)-{1\over \alpha n}\ln E_x^u\left\{e^{\alpha(\sum_{i=0}^{n-1}f(X_i^u)+w_f^{u,\alpha}(X_n^u))}\right\}|+ \nonumber \\
&&|{1\over \alpha n}\ln E_x^u\left\{e^{\alpha(\sum_{i=0}^{n-1}f(X_i^u)+w_f^{u,\alpha}(X_n^u))}\right\}- {1\over \alpha n}\ln E_x^u\left\{e^{\alpha\sum_{i=0}^{n-1}f(X_i^u)}\right\}|+ \nonumber \\
&& |{1\over \alpha n}\ln E_x^u\left\{e^{\alpha\sum_{i=0}^{n-1}f(X_i^u)}\right\}- {1\over n}
E_x^u\left\{\sum_{i=0}^{n-1}f(X_i^u)\right\}|+ \nonumber \\
&& |{1\over n}E_x^u\left\{\sum_{i=0}^{n-1}f(X_i^u)\right\}-\pi^u(f)|= a(\alpha,n)+
b(\alpha,n)+d(\alpha,n)+e(\alpha,n).
\end{eqnarray}
Now notice that: $a(\alpha,n)\leq {\|w_f^{u,\alpha}\|_{sp} \over n}$, $b(\alpha,n)\leq {\|w_f^{u,\alpha}\|_{sp} \over n}$, from Lemma 2 $d(\alpha,n)\leq {n\alpha \over 8}$ and by (UE) $e(\alpha,n) \leq {1 \over n(1-\Delta)} $. Now for a fixed $n$ we let $\alpha$ to $0$ and then let $n\to \infty$ to obtain the claim of Theorem 5.

We shall now consider problem 6.

\noindent
{\bf Theorem 6.} {\it Under (UE), (ME) and (\ref{m1}) for $\alpha_n\to 0$ and ${\cal U}\ni u_n\to u \in {\cal U}$ we have
\begin{equation}\label{fii}
{1\over \alpha_n} I_x^{\alpha_n}(u_n) \to J_x(u)
\end{equation}
as $n\to \infty$.}

\noindent{\bf Proof.} Without loss of generality we may assume that supremum norm of $c$ does not exceed $1$. Them for $f_n(x)=c(x,u_n(x))$, $f(x)=c(x,u(x))$ we have by analogy to (\ref{ossa})

\begin{eqnarray}\label{ossan}
&& |\lambda_x^{u_n,\alpha_n}(f_n)-\pi^u(f)|\leq |\lambda_x^{u_n,\alpha_n}(f_n)- {1\over \alpha_n k}\ln E_x^{u_n}\left\{e^{\alpha_n(\sum_{i=0}^{k-1}f_n(X_i^{u_n})+w_{f_n}^{u_n,\alpha_n}(X_k^{u_n}))}\right\}|+ \nonumber \\
&&|{1\over \alpha_n k}\ln E_x^{u_n}\left\{e^{\alpha_n(\sum_{i=0}^{k-1}f_n(X_i^{u_n})+w_{f_{n}}^{{u_n},\alpha_n}(X_k^{u_n}))}\right\}- {1\over \alpha_n k}\ln E_x^{u_n}\left\{e^{\alpha_n\sum_{i=0}^{k-1}f_n(X_i^{u_n})}\right\}|+ \nonumber \\
&& |{1\over \alpha_n k}\ln E_x^{u_n}\left\{e^{\alpha_n\sum_{i=0}^{k-1}f_n(X_i^u)}\right\}- {1\over k}
E_x^{u_n}\left\{\sum_{i=0}^{k-1}f_n(X_i^{u_n})\right\}|+ \nonumber \\
&& |{1\over k}E_x^{u_n}\left\{\sum_{i=0}^{k-1}f_n(X_i^u)\right\}-\pi^u(f)|=  a(\alpha_n,u_n,k)+
b(\alpha_n,u_n,k)+d(\alpha_n,u_n,k)+e(\alpha_n,u_n,k).
\end{eqnarray}
As in the proof of Theorem 5 we have $a(\alpha_n,u_n,k)\leq {\|w_{f_n}^{u_n,\alpha_n}\|_{sp} \over k}$, $b(\alpha_n,u_n,k)\leq {\|w_{f_n}^{u_n,\alpha_n}\|_{sp} \over k}$ and similarly using Lemma 2  $d(\alpha_n,u_n,k)\leq {k\alpha_n \over 8}$. It remains to estimate $e(\alpha_n,u_n,k)$.  We have
\begin{equation}
|e(\alpha_n,u_n,k)|\leq |{1\over k}E_x^{u_n}\left\{\sum_{i=0}^{k-1}f_n(X_i^u)\right\}-\pi^{u_n}(f_n)|+ |\pi^{u_n}(f_n)-\pi^{u}(f)|=e_1(n)+e_2(n)
\end{equation}
and by Theorem 1 together with Proposition 1 we have that $e_2(n)\to 0$ as $n\to \infty$. From (UE) we have
\begin{equation}
e_1(n)\leq {1\over k} \sum_{i=0}^{k-1}\Delta^i={1\over k(1-\Delta)}.
\end{equation}
Now taking into account from the proof of Theorem 5 that span norms of  $w_{f_n}^{u_n,\alpha_n}$ are bounded by $D$,  for fixed $k$ we let first  $n\to \infty$ and then $k\to \infty$. This way we obtain (\ref{fii}). The proof is completed.

\section{CONCLUSIONS}
In the paper we justify the use of typical approximation procedure for Markov control using piecewise constant functions. Namely we obtain stability (continuity) of functionals assuming that these control functions converge. Using this property we can determine nearly optimal controls using Monte Carlo simulations for finer and finer discretization steps. This fact was known for finite time horizon functionals but there was no studies on that problem for long run functionals.

\end{document}